
\input epsf
\magnification=1200 \overfullrule=0mm \hfuzz=3pt
\catcode`\à=\active  \def à{\`a} \catcode`\é=\active  \def é{\'e}
\catcode`\è=\active  \def è{\`e} \catcode`\ê=\active  \def ê{\^e}
\catcode`\î=\active  \def î{\^\i} \catcode`\ù=\active  \def ù{\`u}
\catcode`\ô=\active  \def ô{\^o} \catcode`\ö=\active  \def
ö{\"o}\catcode`\ä=\active  \def ä{\"a} \catcode`\ï=\active  \def
ï{\"i}\catcode`\ç=\active  \def ç{\c c} \catcode`\û=\active  \def
û{\^u} \font\tenbb=msbm10 \font\sevenbb=msbm5
\newfam\bbfam
\textfont\bbfam=\tenbb \scriptfont\bbfam=\sevenbb
\def\bb{\fam\bbfam\tenbb}
\font\tengoth=eufm10  \newfam\gothfam \textfont\gothfam=\tengoth

\def\qed{\hfill\hbox{${\vcenter{\vbox{ 
    \hrule height 0.4 pt\hbox{\vrule width 0.4 pt height 6pt
    \kern5pt\vrule width 0,4 pt}\hrule height 0.4 pt }}}$}}\line{{\bf Université de Rennes 1}\hfill\hfill}
    \line{{\bf Institut de Mathématiques}\hfill\hfill}
\vskip 20 pt

                       \centerline{{\bf UNE VERSION FEUILLETEE D'UN THEOREME DE BOGOMOLOV}}

\vskip 20 pt
\centerline{FREDERIC TOUZET}
\vskip 10 pt
\hskip 10 pt {\bf Abstract:On compact Kähler manifolds, we classify regular holomorphic foliations of codimension 1 whose canonical bundle is numerically trivial.}
\vskip 20 pt

              \vskip 50 PT

              {\bf 1) Introduction}
              \vskip 5 pt
              1.1 LE THEOREME DE BOGOMOLOV.
              \vskip 5 pt

              Parmi les variétés Kähleriennes compactes $M$, celles dont le fibré canonique est
 numériquement trivial sont caractérisées par le remarquable résultat
              suivant, dont la preuve est due à  F.Bogomolov:
              \vskip 5 pt

              {\bf Théorème 1.1}

              {\it Supposons que $c_1(M)=0$; alors $M$ admet un revêtement fini $\tilde M$ qui est le produit d'un tore par une variété de
              Calabi-Yau.}
              \vskip 10 pt

              Rappelons qu'une variété Kählerienne est dite de Calabi-Yau si elle compacte, simplement connexe et admet une forme volume holomorphe partout non
              nulle.

              Suivant le fameux théorème de Yau ([Ya]) (que nous invoquerons d'ailleurs ultérieurement), on peut munir une telle variété d'une métrique Ricci-plate
              dont l'holonomie induit une décomposition en produit de sous-variétés de Calabi-Yau irréductibles, c'est-à-dire dont le fibré tangent n'est pas
              holomorphiquement décomposable.\vskip 20 pt

              1.2 ANALOGUE FEUILLETE, EXEMPLES ET PRESENTATION DES RESULTATS.
              \vskip 5 pt
              Considérons une variété complexe $M$ munie d'un feuilletage holomorphe régulier de codimension un.
              Ceci correspond à la donnée d'un sous fibré intégrable $\cal F$ du fibré tangent de rang $n-1$ ($n= dim_{\bb C}\ M$).

          Il y a naturellement deux fibrés en droites associés à $\cal F$, le fibré normal et le fibré canonique du feuilletage qui sont respectivement  ${\cal N}_{\cal
          F}={TM\over {\cal F}}$ et $K_{\cal F}=Det\ {\cal F}^*$.
          Par adjonction, on hérite de l'isomorphisme

                             $${{\cal N}_{\cal
          F}}^*\otimes K_{\cal F}=K_M\leqno (1.1)$$
          où $K_M$ désigne comme à l'habitude le fibré canonique de la variété $M$.

          Dans cet article, on s'interesse à la situation suivante: $M$ est Kähler compacte et est munie d'un feuilletage ${\cal F}$ dont la classe canonique est
          numériquement triviale,i.e 
$$c_1(K_{\cal F})=0.$$

          Cette configuration est clairement réalisée par les trois exemples décrits maintenenant:
          \vskip 5 pt
          i) A revêtement fini près, $M$ est le produit d'un tore par une variété de Calabi-Yau et  $\cal F$ est obtenu en tirant en arrière un feuilletage linéaire sur le
          premier facteur par la projection canonique.\vskip 5 pt

          ii)  $M$ est une fibration rationnelle au dessus d'une variété à première classe de Chern nulle (et par suite décrite par le théorème de Bogomolov) et le
          feuilletage ${\cal F}$ est transverse aux fibres.
          \vskip 5 pt
          iii) Le feuilletage ${\cal F}$ est donnée par une fibration en hypersurfaces à première classe de Chern nulle.
          \vskip 5 pt
          Nous affirmons que cette liste est exhaustive:
          \vskip 10 pt
           {\bf Théorème 1.2}

           {\it Soit $M$ une variété kählerienne compacte munie d'un feuilletage holomorphe régulier de codimension un à fibré canonique numériquement trivial; alors ce
           feuilletage est décrit par l'un des exemples (non exclusifs) i),ii) ou iii) ci-dessus.}\vskip 10 pt

{\bf Remarque 1.1}

Si l'on suppose de surcroît que $c_1(M)=0$, on déduit de la
formule d'adjonction (1.1) que ${\cal F}$ est défini par une
forme holomorphe (éventuellement à valeurs dans un fibré plat) qui
est donc identiquement nulle en restriction à toute sous-variété
simplement connexe de $M$. Ceci correspond à la situation décrite
en i). Ce cas facile sera exclu par la suite.
           \vskip 20 pt
           1.3 STRUCTURE DE L'ARTICLE
           \vskip 5 pt
           La première étape consiste à construire, par dualité de Serre, un feuilletage auxiliaire en courbes ( un peu dans l'esprit de [Bo]).

           On montre ensuite que la présence d'un tel feuilletage confère à la variété ambiante des propriétés métriques remarquables (pseudo-effectivité du
           fibré
           canonique). Cette partie s'appuie fortement sur des travaux de Marco Brunella qui étendent dans le cadre Kählerien certains  des résultats
            de semi-positivité établis par Myaoka pour les variétés projectives, et dont une récente version est
             présentée dans [Br]$_2$.

On conclut enfin en déroulant quelques théorèmes classiques de
géométrie riemannienne, en particulier de Yau [Y], Cheeger et
Gromoll [C,G]. \vskip 20 pt {\bf 2) Feuilletage en courbes associé.}
\vskip 10 pt 2.1 PRELIMINAIRES \vskip 5 pt Considérons un
recouvrement ouvert $(U_i)_{i\in I}$ de $M$ tel que sur chaque
ouvert de la famille, les feuilles de ${\cal F}$ soient données par
les niveaux d'une submersion holomorphe

                $$f_i:U_i\rightarrow{\bb C}.$$

Sur chaque intersection $U_i\cap U_j$, les différentielles des
$f_i$ sont liées par la relation

              $$df_i=g_{ij}df_j,\ g_{ij}\in {\cal O}^*(U_i\cap U_j)$$
où le cocycle multiplicatif $g_{ij}$ représente le fibré normal du
feuilletage $N_{\cal F}$.

Le faisceau ${\cal F}_{\infty}$ des germes de $(1,0)$ formes
différentielles lisses tangentes ${\cal F}$ est fin, ce qui assure
l'existence de sections $\omega_i\in{\cal F}_{\infty}(U_i)$
vérifiant sur $U_i\cap U_j$:

            $$\omega_i-\omega_j={dg_{ij}\over g_{ij}}.$$

Puisque les $\omega_i$ sont de la forme
                                           $$g_idf_i$$
où $g_i\in {\cal C}^\infty (U_i)$, on obtient par différentiation
que

       $$\overline{\partial}g_i={g_{ij}}^{-1}\overline{\partial}g_j.$$
On constate ainsi que la collection des $\overline{\partial}g_i$
définit une classe $\alpha$ du groupe de cohomologie
$H_{\overline{\partial}}^{0,1}({N_{\cal F}}^*)$.

En combinant, isomorphisme de Dolbeaut et dualité de Serre, on
obtient que

                 $$Dim_{\bb C}\ H_{\overline{\partial}}^{0,1}({N_{\cal F}}^*)=Dim_{\bb C}\ H_{\overline{\partial}}^{0,n-1}(N_{\cal F}\otimes K_M)$$
et par la formule d'adjonction (1.1)

                $$ H_{\overline{\partial}}^{0,n-1}(N_{{\cal F}}\otimes K_M)=H_{\overline{\partial}}^{0,n-1}(K_{\cal F}).$$

Le fibré $K_{\cal F}$ étant plat, on en déduit par théorie de
Hodge que

                 $$H_{\overline{\partial}}^{0,n-1}(K_{\cal F})\simeq H_{\overline{\partial}}^{n-1,0}(Det\ {\cal F}).$$
\vskip 10 pt 2.2 LE CAS $\alpha=0$. \vskip 5 pt En reprenant les
notations précédentes, il existe sur chaque ouvert $U_i$ une
fonction lisse ${\varphi}_i$ telle que

   $${\overline{\partial}}g_i={\overline{\partial}}{\varphi}_i,\ \ {\varphi}_i={g_{ij}}^{-1}{\varphi}_j\ sur\ U_i\cap {U_j}.$$

On en déduit alors que

               $$\alpha_idf_i-\alpha_jdf_j={dg_{ij}\over g_{ij}}$$
où $\alpha_i=g_i-\varphi_i$ est holomorphe sur $U_i$, en
conséquence de quoi le fibré ${N_{\cal F}}^*$ est muni d'une
connexion {\it holomorphe}.

Puisque la variété ambiante est Kähler, ceci implique la platitude
de  ${N_{\cal F}}^*$ et permet de conclure par la remarque 1.1.
\vskip 10 pt 2.2 LE CAS $\alpha\not=0$. \vskip 5 pt L'isomorphisme
décrit en fin de section 2.1 permet d'exhiber sur $M$ une $n-1$
forme holomorphe $\omega$ non identiquement nulle à valeurs dans
$Det\ {\cal F}$. Localement, une telle forme a pour argument un
multivecteur de degré $n-1$ et de ce point de vue, elle
s'interprète naturellement comme une section holomorphe du fibré
trivial $K_{\cal F}\otimes{Det\ {\cal F}}$. En particulier, la
restriction de $\omega$ à $\cal F$ est soit identiquement nulle,
soit partout non dégénérée. Notons ${\cal C}_\omega$ le
feuilletage en courbe (éventuellement singulier) associé à
$\omega$.

On hérite donc de l'alternative suivante: \vskip 5 pt a) ${\cal
C}_\omega$ est tangent au feuilletage $\cal F$, \vskip 5 pt b)
${\cal C}_\omega$ et $\cal F$ sont partout en position transverse
(en particulier, ${\cal C}_\omega$ est régulier). \vskip 5 pt
C'est cette dernière situation que nous allons considérer dans la
section 3. \vskip 20 pt 3)  METRIQUE INDUITE PAR LE FEUILLETAGE EN
COURBE TRANSVERSE \vskip 10 pt Examinons d'abord le cas où ${\cal
C}_\omega$ admet une feuille rationnelle. Du théorème de stabilité
de Reeb, il résulte que toutes les feuilles sont rationnelles. Le
feuilletage $\cal F$ est donc construit par une suspension de
fibre ${\bb P}^1$ au-dessus d'une variété à première classe de
Chern nécessairement nulle. Il s'agit donc de l'exemple ii) du
théorème 1.2. Notons que l'holonomie de $\cal F$ est ici obtenue
comme représentation de l'extension d'un groupe fini par un groupe
abélien libre dans le groupe de M\oe bius. \vskip 15 pt Si l'on
exclu ce dernier cas de figure, il est établi dans [B,P,T] qu'un
feuilletage en courbe ${\cal C}$ transverse à une distribution en
hyperplan sur une variété Kähler compacte est parabolique ou
hyperbolique: toutes les feuilles sont uniformisées soit
exclusivement par la droite complexe, soit exclusivement par le
disque unité. Il est plus précisément montré dans ({\it loc.cit})
que la distribution transverse est alors intégrable et possède (en
tant que feuilletage) une métrique hermitienne lisse invariante
par holonomie et induisant sur les feuilles de ${\cal C}_\omega$
une métrique naturelle (plate ou de Poincaré suivant le type
conforme de ${\cal C}_\omega$). Dans le contexte étudié, lorsque
${\cal C}={\cal C}_{\omega}$ est parabolique , le fibré ${N_{\cal
F}}^*$ est ainsi muni par transversalité d'une métrique à courbure
nulle et nous sommes à nouveau dans la configuration décrite dans
la remarque 1.1.

Il reste donc à traiter le cas hyperbolique. On obtient de même
que ${N_{\cal F}}^*$ admet une métrique hermitienne dont la forme
de courbure $\Omega$ est une $(1,1)$ forme fermé lisse positive et
invariante par holonomie de ${\cal F}$. \vskip 10 pt {\bf
Proposition 3.1} \vskip 5 pt {\it Supposons le feuilletage en
courbes ${\cal C}_\omega$ hyperbolique.

La variété $M$ peut être alors munie d'une métrique kählerienne
$g$ qui induit sur chaque feuille de ${\cal F}$ une métrique Ricci
plate.

De plus, quitte à changer ${\cal C}_\omega$, $\tilde M$ (le
revêtement universel de $M$) se scinde holomorphiquement et
isométriquement sous la forme ${\tilde N}\times  {\bb D}$ . Cette
décomposition de ${\tilde M}$ est compatible avec celle de $TM$ en
ce sens que $T{\tilde N}$ et $T{\bb D}$ (vus comme sous-fibrés de
$TM$) sont respectivement pull-back de ${\cal F}$ et  ${\cal
C}_\omega$.

Enfin, $\pi_1(M)$ agit diagonalement par isométries sur le produit
${\tilde M}={\tilde N}\times  {\bb D}$.}\vskip 10 pt {\bf Preuve}

Remarquons d'abord que la forme $-\Omega$ est un représentant de
la classe anticanonique $c_1(M)$, eu égard à l'équivalence
numérique de ${N_{\cal F}}^*$ et $K_M$.

D'après le théorème de Yau ([Y]), on en déduit l'existence sur $M$
d'une métrique kählerienne $g$ dont la courbure de Ricci vérifie
l'égalité

           $$Ricc(g)=-\Omega.$$

En chaque point de $M$, le tenseur de Ricci est donc donné par une
forme semi-négative. De plus, par hypothèse de platitude, il
existe sur $M$ une collection de sections locales holomorphes
$\xi_U$ de $Det{\cal F}$ ne s'annulant pas et multiplicativement
liées par un cocycle de module $1$. En particulier, l'ensemble des
$\Delta{{||\xi_U||}^2}$ et ${||\nabla (\xi_U)||}^2$ se recollent
en une fonction globale définie sur $M$ ($\nabla$ désigne ici la
connexion induite par $g$ sur les tenseurs holomorphes de type
$(0,n-1)$). Dans ces conditions, l'inégalité de Böchner
([K,W]p.57)s'énonce dans les termes suivants:

   $$\Delta{{||\xi_U||}^2}\geq {||\nabla (\xi_U)||}^2.$$
Par suite, on obtient que $\nabla (\xi_U)=0$ pour tout ouvert $U$
du recouvrement. En d'autres termes le fibré ${\cal F}$ est
parallèle. Un résultat classique de géométrie kählerienne implique
alors que son orthogonal est {\it holomorphe}. Quitte à
réappliquer les résultats de [B,P,T], on peut conclure que ce
dernier définit un feuilletage {\it hyperbolique}. La
proposition 3.1 résulte alors du classique théorème de
décomposition de De Rham joint au fait que $Ricc(g)$ est
identiquement nulle en restriction à ${\cal F}$.\vskip 10 pt
Précisons maintenant la structure du premier facteur ${\tilde N}$
dans la décomposition du revêtement universel. \vskip 10 pt {\bf
Proposition 3.2} \vskip 5 pt {\it La variété ${\tilde N}$ se
scinde (holomorphiquement et isométriquement) sous la forme ${\bb
C}^k\times V$ où $V$ est une variété de Calabi-Yau.

En outre, $\pi_1(M)$ agit diagonalement par isométries sur le
produit ${\tilde M}={\bb C}^k\times V\times  {\bb D}.$}\vskip 10
pt

La preuve repose sur le résultat suivant de Cheeger et Gromoll.
\vskip 10 pt {\bf Théorème 3.1} ([C,G])

{\it Soit $M$ une variété lisse simplement connexe admettant une
métrique riemannienne complète de courbure de Ricci positive ou
nulle; alors $M$ se décompose isométriquement sous la forme ${\bb
R}^k\times V$ telle que $V$ ne contient pas de droite géodésique.}
\vskip 5 pt Rappelons qu'une droite géodésique est une géodésique

$$ \gamma:]-\infty, +\infty[\rightarrow M$$
telle qu'à chaque instant $t,t'\in{\bb R}$,$\gamma_{|[t,t']}$ soit
minimale. \vskip 5 pt {\bf Preuve de la proposition 3.2}

Puisque la décomposition de De Rham d'une variété kählerienne
coïncide avec celle de la variété réelle sous-jascente ([K,N] th
8.1, p 172], on obtient que ${\tilde N}$ se scinde
holomorphiquement et isométriquement sous la forme ${\bb
C}^l\times V$.

Par construction, le groupe fondamental de $M$ agit diagonalement
par isométries sur le produit ${\tilde M}={\bb C}^l\times
V\times{\bb D}$. Il reste à prouver que $V$ est compacte;
supposons par l'absurde que ce ne soit pas le cas. La contradiction
recherchée résulte alors de l'argument suivant  dû à Cheeger et
Gromoll ({\it loc.cit}). Soit $K$ un domaine fondamental pour
l'action de $\pi_1(M)$ sur ${\tilde M}$. C'est un compact et sa
projection $\rho (K)$ sur $V$ est donc un compact dont l'orbite
par $\rho (\pi_1(M))$ est $V$ toute entière. Puisque $V$ est
supposée non compacte, il existe en un point $p\in \rho (K)$ un
rayon géodésique $\gamma:[0,\infty[\rightarrow V$(i.e une
demi-droite géodésique) tel que $\gamma (0)=p\in \rho (K)$. Soit
$g_n$ une suite d'isométries de $\rho (\pi_1(M))$ telle que
$g_n(\gamma(n))=p_n\in  \rho (K)$. La géodésique $\gamma_n$ de $V$
définie par $\gamma_n(0)=p_n$ et
${\gamma_n}^{'}(0)=dg_n({\gamma}^{'}(n))$ est donc un rayon
géodésique en restriction à $[-n, \infty[$. Par compacité, on peut
 extraire une sous-suite $n_i$ de sorte que $p_{n_i}$ et
${\gamma_{n_i}}^{'}(0)$ convergent respectivement vers $p_0\in
\rho (K)$ et $v_0\in T_{p_0}^1 (V)$. Par construction, la
géodésique passant en $p_0$ à la vitesse $v$ est une droite
géodésique.\vskip 20 pt 4) PREUVE DU THEOREME 1.2.\vskip 10 pt

Au vu des sections précédentes, il reste à traiter les cas  \vskip
5 pt

(1) ${\cal F}$ admet un feuilletage en courbes hyperbolique transverse

(2) ${\cal F}$ contient un feuilletage en courbe (éventuellement
singulier) défini globalement par $\omega\in\Omega^{n-1}(Det\
{\cal F})$.

Examinons d'abord le point (1). En reprenant les notations de la
section 3, il existe alors une variété de Calabi-Yau $V$ telle que

                          $${\tilde M}={\bb C}^l\times V\times {\bb D}$$
avec une action diagonale et isométrique de $\pi_1(M)$.

\vskip 5 pt

{\bf lemme 4.1} {\it L'action de $\pi_1(M)$ sur ${\bb D}$ est
discrète.} \vskip 10 pt
 Puisque le groupe des isométries de $V$ est fini ([B]), on est ramené à considérer un groupe $H$ dont l'action est  libre,
 discrète, diagonale, cocompacte et isométrique sur le produit $ {\bb C}^l\times {\bb D}$. Il reste à voir que cette action reste
  discrète sur le second facteur, auquel cas chaque feuille de $\cal F$ sera fermée, situation décrite par l'exemple iii) du théorème 1.2.

Ceci va résulter du fait général suivant: \vskip 10 pt {\bf
Théorème 4.1 (Auslander)} (cf.[R] p 149)

{\it Soit $G$ un groupe de Lie et $R$ un sous-groupe connexe,
résoluble et distingué. Soit $\pi:G\rightarrow G/R$ la projection
canonique. Soit $H$ un sous-groupe fermé de $G$ dont la composante
neutre $H^0$ est résoluble. Alors la composante neutre de
$\overline{\pi (H)}$ est résoluble.} \vskip 5 pt {\bf Preuve du
lemme 4.1} Il suffit manifestement de montrer que $H$, vu comme
réseau de ${\bb C}^l\times U(l)\times Sl(2,{\bb R})$ se projette
sur un {\it réseau} de $Sl(2,{\bb R})$. On note $H_1$ l'image de
$H$ par cette projection; remarquons que $H_1$ n'est pas conjugué
à un sous-groupe triangulaire de $Sl(2,{\bb R})$ car ${{\bb
C}^l\times {\bb D}}/H$ est de première classe de Chern non nulle. Par
conséquent, $H_1$ est, soit un réseau, soit dense dans $Sl(2,{\bb
R})$. Supposons par l'absurde que $\overline{H_1}= Sl(2,{\bb R})$.
Notons $\pi$ et ${\pi_1}$ les projections canoniques respectives
de  ${\bb C}^l\times U(l)\times Sl(2,{\bb R})$ sur $U(l)\times
Sl(2,{\bb R})$ et de $U(l)\times Sl(2,{\bb R})$ sur $Sl(2,{\bb
R})$. D'après le théorème d'Auslander, la composante neutre
$\Gamma$ de $\overline{\pi (H)}$ est résoluble; par suite,
$\pi_1(\pi (\Gamma))$ est un sous-groupe (non nécessairement
fermé) résoluble de $Sl(2,{\bb R})$. Par hypothèse, l'adhérence de
$\pi_1(H)$ est précisément $Sl(2,{\bb R})$; ceci entraîne alors
que $\pi (H)$ admet une infinité de composantes connexes
${(\Gamma_n)}_{n\in{\bb N}}$ telles que $\bigcup_{n\in{\bb
N}}\pi_1(\Gamma_n)$ soit dense dans $Sl(2,{\bb R})$. Il en résulte
que $\pi_1(\pi (\overline{H}))=Sl(2,{\bb R})$ par compacité de
$U(l)$, ce qui est impossible, vu que $Sl(2,{\bb R})$ ne peut
s'exprimer comme union dénombrable de translatés d'un sous-groupe
résoluble.

\vskip 10 pt Il reste à analyser le point (2). Selon la
terminologie utilisée par Brunella ([Br]$_2$), le feuilletage
${\cal C}_\omega$ défini par $\omega$ n'est pas un feuilletage en
courbes rationnelles. En effet , dans le cas contraire, le degré
de $K_{\cal F}$ sur une droite générique de
${\cal C}_\omega$ serait strictement négatif.

Dans ces conditions , Brunella ({\it loc.cit}) a établi que le
fibré cotangent de ${\cal C}_\omega$ est pseudo-effectif. Par
ailleurs, la $n-1$ forme $\omega$ est holomorphe à valeurs dans un
fibré plat; par adjonction, on obtient alors que $K_M$ et donc
${N_{\cal F}}^*$ sont pseudo-effectifs.

Le lemme qui suit a été montré par Brunella ([Br]$_1$) dans le
cadre des feuilletages sur les surfaces complexes compactes et
s'adapte sans difficultés lorsque la variété ambiante est
kählerienne. Nous en donnerons la démonstration par commodité pour
le lecteur (voir aussi [B,P,T]).

\vskip 10 pt {\bf lemme 4.2} {\it Soit $\cal F$ un feuilletage
holomorphe de codimension 1 sur $M$ Kähler compacte. Supposons que
${N_{\cal F}}^*$ soit pseudo-effectif; alors il existe un courant
positif fermé de bidegré $(1,1)$ invariant par $\cal F$.} \vskip
10 pt {\bf Preuve:} Par hypothèse, le fibré ${N_{\cal F}}^*$ admet
une métrique singulière dont les poids locaux sont
plurisousharmoniques. La forme de courbure de cette métrique est
donc donnée par un courant positif $T$ localement défini par
$T={\sqrt{-1}\over 2\pi}\partial\overline{\partial}F$, où $F$ est
un poids local de la métrique. On peut choisir ces poids locaux de
façon à obtenir sur $M$ une forme $\eta$ (à priori non lisse)
telle que, localement:

$$ \eta=\sqrt{-1}e^{2F}df\wedge d\overline{f}.$$
où $\{df=0\}$ est une équation locale de ${\cal F}$. Soit $\theta$
une forme de Kähler; on déduit du théorème de Stokes (pour les
courants) que

   $$\int_M \partial\overline{\partial}\eta\wedge{\theta}^{n-2}=0$$

D'autre part, le courant $\partial\overline{\partial}\eta$ est
positif (l'exponentielle d'une fonction $psh$ est $psh$). La
nullité de l'intégrale précédente implique alors que $e^F$ est
pluriharmonique et donc constante {\it dans les feuilles}. \vskip
10 pt

L'existence d'un tel courant invariant par holonomie  confère au
feuilletage $\cal F$ des propriétés qualitativement semblables à
celles des feuilletages riemanniens. Plus précisément, il est
montré dans [Br]$_1$ que ${\cal F}$ entre dans la catégorie des
feuilletages quasi-uniformément isométriques étudiés par Kellum
([K]). Suivant Brunella ({\it loc.cit}) , on peut alors affirmer
que $M$ est union disjointe de minimaux de $\cal F$, chaque
minimal étant soit une feuille compacte, soit une sous-variété
réelle topologique de codimension 1, soit $M$.

Si le feuilletage est minimal, on obtient par transitivité(voir
[Br]$_1$ pour les détails) que $\cal F$ admet une métrique lisse
transverse invariante par holonomie et de courbure négative
constante. De même que dans la preuve de la proposition 3.1, on
peut alors invoquer le théorème de Yau pour conclure à l'existence
d'un feuilletage {\it holomorphes} en courbes {\it transverse}.

Si l'on élimine le cas des fibrations ( iii)du théorème 1.2), il
reste à étudier les feuilletages comportant un minimal de
codimension 1 réelle. Ceux ci sont décrits comme suit ([Br]$_1$):

a) l'adhérence de chaque feuille est de codimension 1 réelle.

b) il y a exactement deux feuilles compactes et les autres
feuilles ont des adhérences de codimension 1. \vskip 10 pt Rappelons
maintenant un des principaux résultats de Kellum: dans le cadre
quasi-isométrique, l'adhérence du pseudo-groupe d'holonomie est un
pseudo groupe de Lie. En particulier, on hérite d'un faisceau
localement constant de germes de champs de vecteurs transverses à
$\cal F$. Par un argument qu'on peut trouver par exemple chez
[Lo,Re], l'ensemble de ces germes est en chaque point une algèbre
de Lie (réelle) de champs de vecteurs holomorphes. Cette algèbre
est visiblement de dimension $1$ dans les cas a) et b) ci-dessus.
Par dualité, le feuilletage $\cal F$ est donné par une forme
logarithmique (éventuellement à valeurs dans un fibré plat) sans
diviseurs de zéros et dont l'ensemble des  pôles est précisément
formé de la réunion des feuilles compactes. On élimine ainsi le
cas b), incompatible avec la pseudo-effectivité du fibré conormal.
Dans la situation a),  ${N_{\cal F}}^*$ est numériquement trivial
(exemple i) du théorème 1.2).

\vskip 50 pt

{\bf BIBLIOGRAPHIE}
\vskip 20 pt

\parindent=-1cm\leftskip=-\parindent
[B] A.Beauville; {\it Variétés Kähleriennes dont la première
classe de Chern est nulle}, J.Differential Geom.18 (1983), 4,
755-782.\par [Bo] F.A Bogomolov; {\it Kähler manifolds with
trivial canonical class}, Math.USSR Izvestija, Vol.8 (1974), $N^o$
1.\par [B,P,T] M. Brunella, J.V. Pereira, F.Touzet; {\it Kähler
manifolds with split tangent bundle}, Bulletin de la
SMF {\bf 134}, fasc. 2 (2006).\par [Br]$_1$ M. Brunella; {\it Feuilletages holomorphes sur
les surfaces complexes compactes}, Ann.scient.Ec.Norm.Sup., $4^e$
série,t.30,1997, p.569-594.\par [Br]$_2$ M.Brunella; {\it A
positivity property for foliations on compact Kähler manifolds},
Intern.J.Math. 17 (2006), 1, p. 35-43.\par [C,G]
J.Cheeger,D.Gromoll; {\it The splitting theorem for manifolds of
nonnegative Ricci curvature}, J.differential geometry 6 (1971/72),
119-128.\par [K] M.Kellum;{\it Uniformly quasi-isometric
foliations} Erg.th.Dyn.Sys,Vol.13 (1993) P.101-122.\par [K,N]
S.Kobayashi, K.Nomizu; {\it Foundations of differential
geometry.Vol.II.} Wiley \& sons, Inc.,New York,1996.\par [K,W]
S.Kobayashi,H.Wu;{\it Complex differential Geometry}, DMV Seminar
Band 3, Birkhäuser,(1987).\par [L,R] F.Loray, J.C.Rebelo;{\it
Minimal, rigid foliations by curves on} ${\bb C}{\bb P}(n)$,
Journ.Euro.Maths.soc.5 (2003),p.147-201.\par [R] M.S Raghunatan;
{\it Discrete subgroups of lie groups}, Erg.der.Math, 68,
Springer, (1972).\par [Y] S.T. Yau;{\it On the Ricci curvature of
a compact Kähler manifold and the complex Monge-Ampère
equation},I, Comm.Pure Appl.Math.31 (1978),339-411.

\end